\documentclass{amsproc}
\usepackage{graphicx}
\usepackage{hyperref}
\usepackage{appendix}
\usepackage{amscd}
\usepackage{amsmath}
\usepackage{epsfig}
\usepackage{amsfonts}
\usepackage{amssymb}
\usepackage{hyperref}
\usepackage{pdfpages}
\usepackage[thinc]{esdiff}
\usepackage{accents}
\usepackage{nicefrac}

\theoremstyle{definition}

\theoremstyle{remark}

\numberwithin{equation}{section}
\allowdisplaybreaks

\input{tcilatex}

\begin{document}
\title[On Dick Askey]{My Warm, Randomly Recorded, Recollections \\
of Professor Richard Askey\/}
\author{Sergei K. Suslov}
\address{School of Mathematical and Statistical Sciences, Arizona State
University, P.~O. Box 871804, Tempe, AZ 85287-1804, U.S.A.}
\email{sergei@asu.edu}
\date{today}

\begin{abstract}
These are my memories of moments with Dick and Liz Askey in Russia, Wisconsin, Arizona, and abroad. 
Dedicated to the Askey family, these recollections span over 40 years and encompass many dramatic changes in the world. Due to this, it is challenging to entirely separate my personal thoughts and feelings from the factual historical account.
\end{abstract}

\maketitle


Dick Askey's main scientific legacy lies in the theory of orthogonal polynomials and special functions. 
{\it{De facto\/}}, he renewed worldwide interest in this classical area of analysis in the second half of the last century \cite{Askey75} and since then had been one of its international leaders for several decades \cite{An:As:RoySF}. The introduction of the Askey--Wilson polynomials made him a classic author \cite{AW85} (see also \cite{At:SusAW}). Dick was a natural ``connector" among experts, young people, and students from different countries. 
His true legacy will be much better understood later \cite{Liber, RAleg}.

As many others in this area, I was very lucky to consider Askey as a
friend/sensei since the late 1980s.
In a way, he, along with Mourad Ismail, Mizan Rahman, and Joaqu\'{i}n Bustoz, 
drastically changed my life -- from being a junior researcher 
in theoretical and mathematical physics at the Kurchatov Institute of Atomic Energy 
in Moscow to securing a professorship in mathematics at Arizona State University. 
This transformation was only possible during the disintegration of the Soviet Union 
and the opening of the ``Iron Curtain", the political and ideological barrier 
that separated the USSR and its satellite states in Eastern Europe 
from the rest of the world during the Cold War%
\footnote{ I may also share Dick's opinion about the end of the Cold War: he said once that 
it (~the``Iron Curtain") was so solid, made forever, but then suddenly disappeared just in a few months.  
}
.
Those are my recollections of some episodes over the years that in my opinion reveal 
Askey's personality not only as a great mathematician but also 
as a great human being (and a true friend). 
Thank you very much, Dick, for all you have done for me and 
for our numerous colleagues and friends all over the world!

\section{In Moscow, at the Kurchatov and Steklov Institutes\/}

I met Professor Askey for the first time during his very first visit to Moscow through 
an academic exchange program in the late 1980s \cite{Tom}. That day, he was picked up by Arnold Fedorovich Nikiforov and his daughter Lyba, along with Professor Vasilii Borisovich Uvarov, from the academic hotel.%
{\footnote{
According to the travel diaries by Liz Askey, this meeting `with a trio of physicists' took place on Monday, September 14th, 1987 \cite{Tom}.
She also writes: ``One of the physicists who came by this morning brought me a bouquet of flowers
from his garden.'' It must have been Nikiforov, because at this time of year he lived at his ``dacha'' near Moscow. 
}}
She had excellent English and was our translator. 
Indeed, Dick was the first foreigner I ever met. At the Moscow Institute of Physics and Technology \cite{FizTex},
we had extensive English classes in small groups of students, but we were mainly learning how to read and write -- obviously, mostly technical literature -- and our spoken English was only practiced among us, fellow students, and with our instructors, who were also educated in the Soviet Union.
Therefore, for me,
who studied German in the high school,
it was hard to understand Dick during our first meeting. 
For instance, it took me some time to figure out Dick's correct pronunciation 
of the word ``measure" because I had never heard it before in conversation. 
But our ``communication in Latin" was OK; Dick immediately recognized 
my ``\textit{mutantis mutandis}", the well-known expression we had already used 
in writing in our article on the continuous Hahn (actually Chebyshev, 1875) polynomials \cite{At:SusCH} 
(see also \cite{Askey85}%
%
%
\footnote{%
In a reprint of this article that Dick gave to us in Moscow, he penciled in ``I should have done this if I had been smart enough!" to the best of my recollection.}%
).
Much later, Dick mentioned to me several times that we always understood each other.

\begin{figure}[h!]
\includegraphics[width=0.765\linewidth]{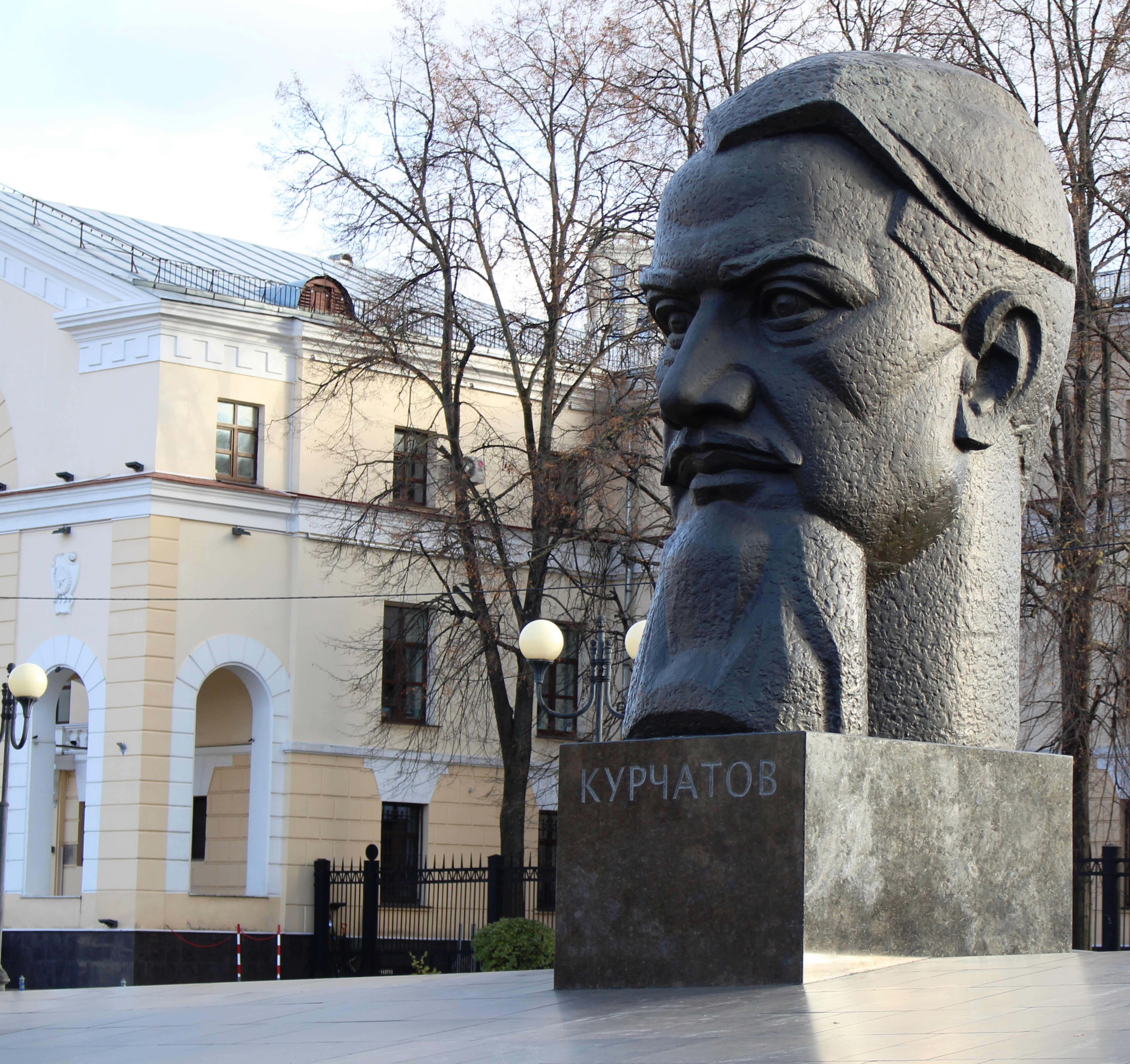} 
\caption{ Russian Research Center ``Kurchatov Institute", Moscow. }
\end{figure}
Dick and the rest of the company came to Kurchatov Square in Moscow by trolleybus from the subway, 
which stopped next to a bust of Kurchatov’s head in front of the institute's main entrance (Figure~1, \cite{Kurchatov}). He introduced himself by saying, ``Askey", and shook my hand. As a foreigner, at that time, Dick was not allowed to enter the Kurchatov Institute, but there was a special house just across the street -- a~two-story international office building -- where we could meet. 
Anyone who worked for the institute had to follow this procedure. 
After entering a narrow meeting room on the second floor, Dick immediately remarked that ``the table is too long 
and the blackboard is too small". Over the years, this became a typical reaction from Dick to ``real-life circumstances". 

I was not heavily involved in his conversation with Nikiforov and Uvarov. 
At that time, we worked on our approach to the theory of classical orthogonal polynomials of a discrete variable, through the difference equation they satisfy \cite{NSU}, and I tried to connect the ``Russian and American ways".
I was mostly trying to understand what Dick was saying and, alternatively, what he was not saying. 
It became clear to me later that, at that time in Moscow, he never mentioned the biorthogonal rational functions \cite{Wilson91}.

In addition, there was a meeting at the Steklov Institute in a small office of the Academician Sergey Mikhailovich Nikol'skii, who was my advanced calculus instructor at the Moscow Institute of Physics and Technology. 
Most of our math professors were from the Steklov Institute \cite{SteklovMoscow}, 
which helped me a lot later in my career dealing with mathematical journals. 
After an informal seminar with Dick as the main speaker, the Academician Andrei Aleksandrovich Gonchar, later a vice president of the Russian Academy of Sciences, invited this small group of 
interested people to lunch at his apartment, which was very close to the institute. 
{\it{De facto\/}}, because of Dick, I was invited to the ``private company" of the future vice president 
of the Academia!

In the evening, Dick also gave a talk at the famous Israel Moiseevich Gel'fand seminar at Moscow State University, but I could not get through the security there, in the main university building, and unfortunately had to miss this talk. Gel'fand and his coauthors were developing their general theory of hypergeometric functions at that time; see, for example, \cite{GelGr} and the references therein. It is worth noting that from the very beginning of his first trip to Moscow, Professor Askey said at the Academy of Sciences that he would like to meet the Academician Gel'fand and Suslov during his visit, which was a great surprise and honor to me. I was just a humble young researcher at the ``highly classified" Kurchatov Institute, who never expected to meet a foreigner! Personally, I feel that Dick somehow liked me from our first meeting, although there was not much conversation, and those feelings were mutual for the rest of his life.

\section{In Florida and Sicily\/}

My first visit abroad, to the USA and Canada, occurred thanks to Dick Askey and Mourad Ismail, who, together with Ed Saff, organized an International Conference on Approximation Theory in Tampa, Florida, in March 1990. This was because of a new policy by Gorbachev, which made it possible for those working at the Kurchatov Institute to attend conferences abroad. I had been invited before to East Germany and Bulgaria, but the USSR Atomic Energy Committee never gave me permission to visit those socialist countries. As a result of ``perestroika", I visited the United States and Canada before ever going to any socialist country -- probably the first such case at the Kurchatov Institute in the 1990s!

Our small group of participants from the Soviet Union, mainly from the Steklov Institute, was surprised by the American hospitality and the exotic beauty of Florida! Dick and Liz were very friendly; they took two of us, 
namely, Natig Atakishiyev and me, 
for a lovely lunch one of those days. In a plaza, while walking to the restaurant from the parking lot, Dick suddenly stopped at the door of an abandoned office and said that the US also had economic problems. It was his usual way to support people, {\it{de facto\/}}, saying that things should go better in the USSR, or, at least, I understood his comment this way.

Various people told me later that Dick was a ``straight shooter" who never hesitated to tell the truth, if something wasn't right, but I never saw him that way! For me, he was always supportive, encouraging, and optimistic! (Possibly because I never said something really stupid in his presence?%
\footnote{%
On a contrary, I once mentioned to Dick -- who needs a million dollars --
meaning that the life is good in the United States if one has a decent job.
But he immediately replied that he would like to have 10 million dollars in order
to give 5 million dollars to the Chudnovsky brothers! \cite{Chud}\/} 
) 
Dick also took us on a lovely trip to one of Florida's parks, with wooden bridges and benches on the trails, turtles in the lake, exotic birds, palm trees, etc. 
Mourad and Thanaa had an unforgettable reception at their house in Tampa.
After the meeting, Natig Atakishiyev and I went to Ottawa by a kind invitation of Mizan Rahman \cite{Rahman}.
These visits created a solid foundation of our fruitful collaboration for many years to come.
%

%
Later that spring, we met once again in Erice, Sicily, a beautiful place on the top of a mountain with very narrow, wet, or even icy streets due to the elevation. I remember how amazed Dick was by a skillful driver who brought him and Liz from the local airport, driving very fast just a few centimeters away from solid walls! Like the Academician Lev Davydovich Landau, my scientific grandfather (Landau--Smorodinsky--Suslov \cite{Smor}, \cite[Contribution \#~71]{Liber}, 
Figure~2), Dick really liked people who did their job very well! (Paul Nevai, one of the organizers of the meeting, the Third International Symposium on Orthogonal Polynomials and Their Applications, Erice, Italy, June 1990, with whom I met for the first time in Florida as well, should remember this episode?)
A light, gray conference bag from Erice became Dick's `travel companion' for many years to come.
By recalling Erd\H{o}s number as a level of collaboration in mathematics,
Askey's Erd\H{o}s number is 2, and therefore my Erd\H{o}s number is 3 (Erd\H{o}s--Boas--Askey--Suslov, 
one of the paths) which is due to Dick. 
He always called us physicists \cite{Tom}!
Thus Askey's own ``Landau number" is also 3 (Landau--Smorodinski--Suslov--Askey) as a result of our collaboration \cite{As:Su1, As:Su2}.
Actually this is an example that mathematics and physics are deeply interviewed and rely on each other in order to progress.
\begin{figure}[h!]
\includegraphics[width=0.657\linewidth]{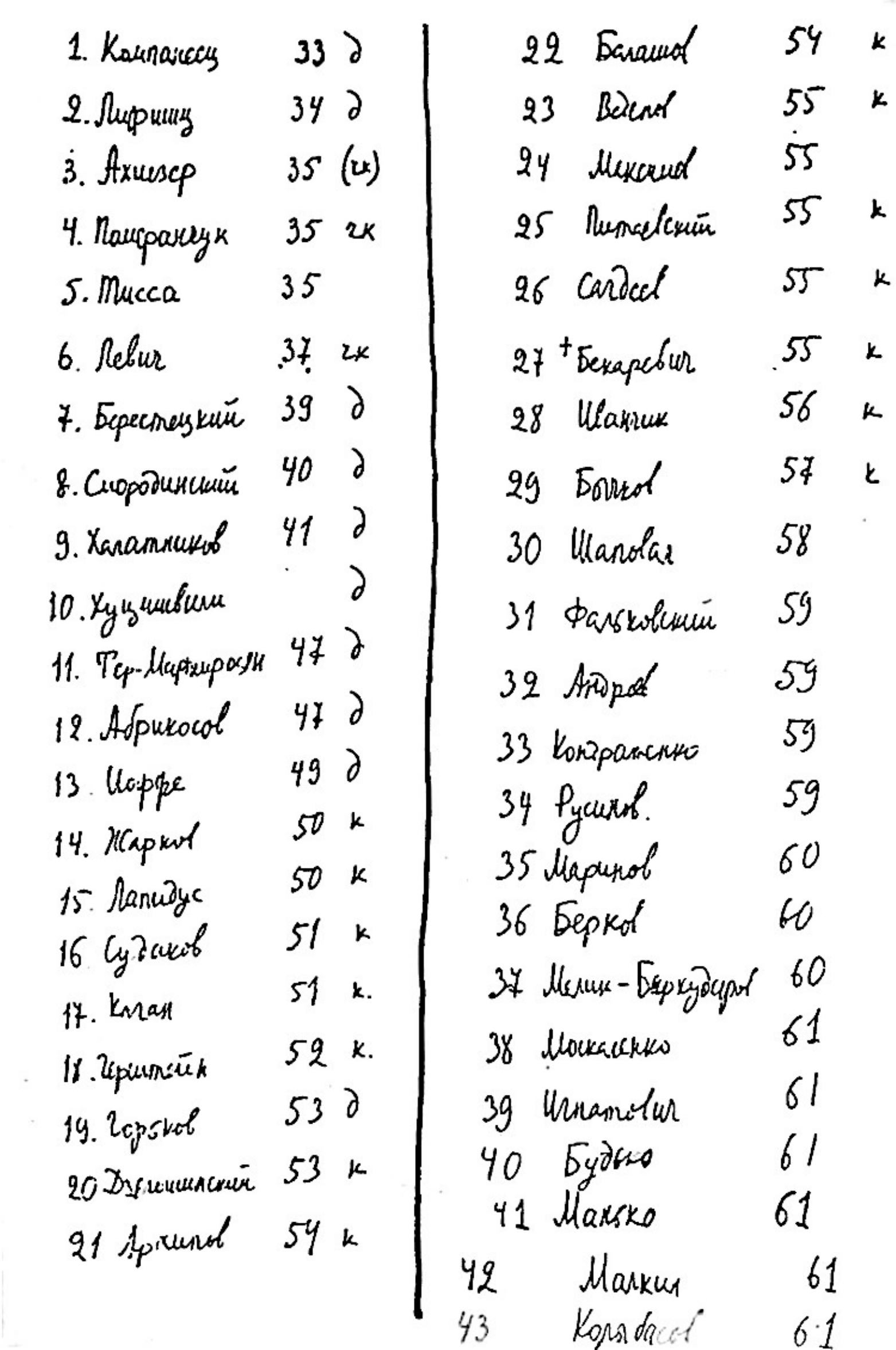} 
\caption{ A copy of Landau's Theoretical Minimum list; among passing the
exams, \#~8 is Professor Yakov Abramovich Smorodinski, in~1940, Doctor of Science; my master thesis and Ph.~D.
advisor; see \cite{Lan:Lif} for a true story of Landau's exams.}
\end{figure}

\section{In Leningrad and Moscow\/}

The meeting in Florida continued the following year in Leningrad, now known as Saint Petersburg. The local branch of the Steklov Institute and the Euler International Mathematical Institute \cite{SteklovStP} organized a two-week program, Methods of Approximation Theory in Complex Analysis and Mathematical Physics, in May 1991, in a newly renovated building of this institute on the Fontanka River. 
The activity of the Euler Institute consists in the organization and running of scientific programs, conferences and workshops on mathematical problems of fundamental importance and related topics.
Academician Ludwig Faddeev was the director of this branch of the Steklov Institute and the founder of the Euler Institute \cite{Faddeev}.
He attended some lectures and recorded some episodes on video camera, 
but it might had been at another meeting, namely, on a workshop dedicated to development of quantum groups. 
Lawrence Biedenharn was one of the distinguished visitors of  the Euler Institute from abroad, along with Dick Askey, young Masatoshi Noumi, and Mourad Ismail.
Somehow, Dick arrived only for the second week of this two-week meeting. However, his flight was scheduled to land at Moscow International Airport, and he needed to catch a connection from Moscow to Leningrad. In the West, this might have been easy, but in those days in Moscow, it was a real problem for a foreigner. Dick could easily get lost in the big city, especially with a short interval between flights, and no one around could speak English. This is why I had to take an overnight train from Leningrad to Moscow to pick up Dick at the airport and then go back.
Luckily, his international flight was on time and there was no problem with customs. As far as I remember, 
I told him ``Follow Me" -- the title of an international airport journal in Italy. Our local bus journey went smoothly, but it was probably crowded by Western standards; there were plenty of similar Hungarian made Icarus buses later in Canada though. Dick had a big yellow bottle with water, which he wanted to take aboard. On the way to the second local terminal, he asked me how to say ``water" in Russian because over there, at the local airport security, no one was supposed to speak English!
Luckily all went well. 
His second flight to Leningrad was smooth, and we had a great time in St.~Petersburg including a wonderful tour of the Peterhof Palace, built by some of the Europe's best architects of the 18th century \cite{Petergof}.
After the meeting, Mourad and Dick came back to Moscow for a few extra days
{\footnote{%
One of the vivid recollections was the visit with our American colleagues to
Treyakov Gallery, where ``we had a moment" with Dick in front of the famous
home portrait of Nicholas~II, by Valentin Serov, at that time just appeared on
display (Figure~3). In my opinion, and Dick agreed: Looking at the last
Russian emperor's eyes, one may say that the man knew his fate.}} .

\begin{figure}[h!]
\includegraphics[width=0.69%
\linewidth]{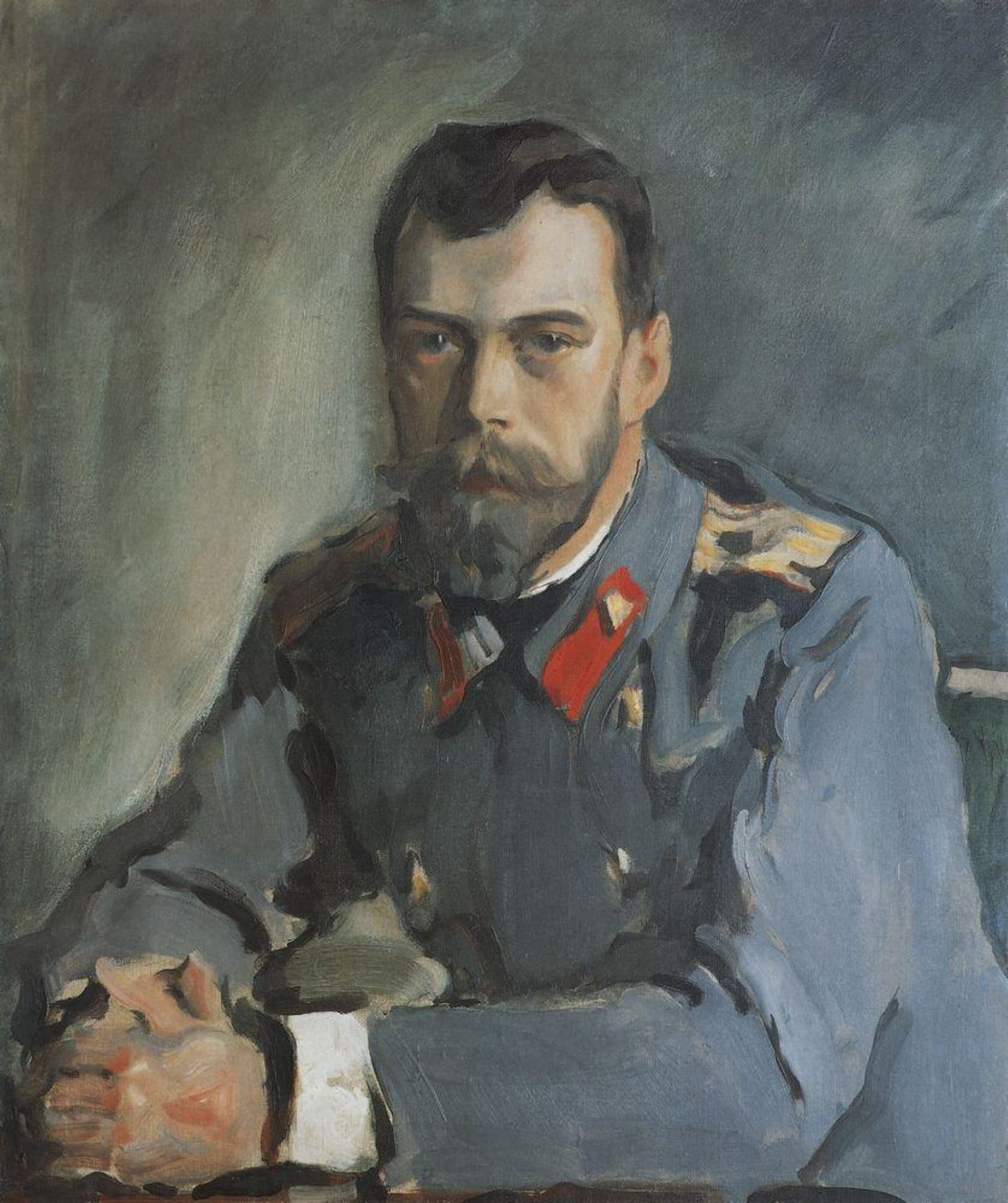} 
\caption{Portrait of the last Russian emperor Nicholas~II, by Valentin
Serov, 1900; Tretyakov Gallery, Moscow. }
\end{figure}

Many years later in Arizona, when Dick was in Tucson for a meeting on
mathematical education, I came once again from Phoenix to take him to the
airport over there, just to enjoy his company on a 20-minute ride from his
hotel and during a lunch at the airport.

\section{In Madison, with Dick and Liz\/}

I had the privilege of staying at Dick and Liz's house for about a week, allowing me to observe their professional and personal lives from within the context of ``regular everyday life routine". 
The flight from Canada was on time, actually a bit early, but they picked me up at the airport without any delay.
Passing by the Oscar Mayer factory, Askey said that it's the only production line around -- there are no other factories in Madison.
My first-ever Thanksgiving was at Dick's house in November 1993, shared with the entire Askey family. Dick woke up early that day and did all the cooking: turkey, soup, stuffing, etc. The remaining meat went into our lunch sandwiches the next week, which Dick prepared himself for both of us every morning. Professor Walter Rudin's office was next to Dick's at the university, and I received his book on complex analysis as one of my ``trophies".
Many faculty members met for a lunch in the math lounge with a beautiful view of the lake and the city. 

The amazing library at the math department had a great collection of rare books on orthogonal polynomials and special functions. Shaun Cooper, Dick's Ph.~D. student at that time, who had just received a teaching award for excellence, immediately solved a characterization problem for our nonuniform lattices, which I later included in my book as an exercise \cite{SusFS} (with a proper credit, of course). 
Dick was very happy for Shaun because, in his own words, ``those awards are not given easily". 
In addition, Shaun was completing an article on a characterization of elliptic functions. 
In the university bookstore, on the shelf, I found the book \cite{AtomicBomb}, which contained 
the highly classified, at that time, history of the Kurchatov Institute.
I had a chance to give a talk at Dick's seminar on the ``tree method", an elegant graphical approach to the solution of the multi-dimensional Laplace equation, invented by Vilenkin, Kuznetsov, and Smorodinsky; see \cite{NSU} and the references therein. Dick carefully read and edited my hand-written lecture notes in the living room of his house and then distributed those among the students; he was working with one of his female calculus students on  evaluation of the volume of $n$-dimensional sphere as a honors project. I am using these notes for my quantum mechanics class at Arizona State University till now! 
Shaun  \cite[Contribution \#~19]{Liber} provides more details on Dick's teaching in Madison.
Dick was running his household as a ``well-oiled machine", energetically moving from the basement to the kitchen, to the second floor and back, everyday starting early in the morning, before walking to the university. 
Only with this in mind, one can truly understand the well-known expression that
``the best nap is during somebody's talk" that he practiced a lot!
All household supplies came from the local neighborhood store located on the corner. Dick liked to go there and he introduced me to the owner. 
Relationships with their neighbors were very friendly. One day,  about 4:30 in the morning, his next door neighbor phoned in and asked if Dick realizes that his car is parked on the street in front of the house. Dick yelled `oh' and immediately went outside in order to re-parked the car in the drive way. It was a ``cleaning day" on their side of the street and he could of get a fine otherwise!
In Madison, I was also invited to another private event: Liz and Dick took me to a gathering that was unusual for me, held in the basement of their church (Trinity United Methodist Church \cite{Trinity}?). Everyone brought some pasta, salads, etc. To my surprise, Dick was an elected leader of this local religious community, delivering an opening prayer on this dinner. To me, it was an example of the total respect people had for him, those who were far away from our mathematics community, namely, their neighbors in Madison. Everyone was friendly to me, a nonreligious person by the Soviet education. Someone told me that their church had collected donations and food for people in the former USSR, which was in bad economic shape at that time.
I should also share my humble observations on the relationship between Dick and Liz: he loved her deeply and always cared for her. When I visited their home, every time Dick returned from work, he would tell Liz everything that happened during his day, and she always listened. She was deeply involved in all aspects of his life! There was a strong bond between the two of them; they had always been a family! Liz worked at a children's bookstore, and visiting bookstores in Madison was one of their family's traditions. Mourad also told me that Liz was a judge of the perfection of English in their articles under preparation; see also \cite{Tom} for her own excellent writing. Later, when Dick became interested in mathematics education, her phrase ``I married you for better or for worse, but not for the mathematical education" became legendary in our society. Dick was very proud of his grandson, Ben, who was a boy at that time; they both loved all their grandchildren very much.
\begin{figure}[h!]
\includegraphics[width=0.897\linewidth]{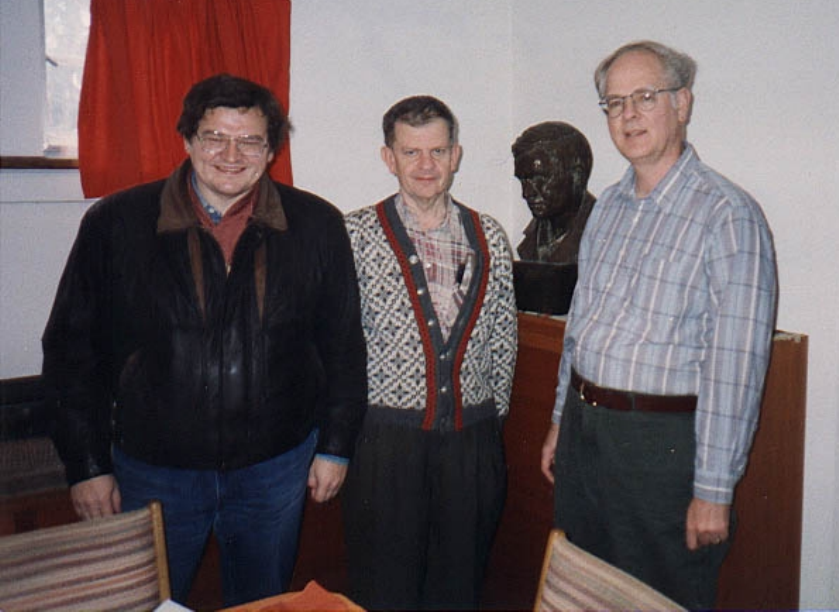} 
\caption{ Sergei Suslov, Dick Askey, Ramanujan's bust, and George Gasper in
Dick Askey's house in November 1993. }
\end{figure}
In the house, they put me in their daughter Suzanne's unoccupied room, which also contained Dick's personal computer, on the second floor. Every morning, Dick was typing with both hands, without looking at the keyboard, answering numerous emails before going to work -- he had the skills of a professional typist! There was also a box full of articles from different mathematical journals, waiting for reviews, and every evening Dick refereed one of them with more articles coming in... without any hope that this box would ever be empty -- there is only so much you can do!%
{\footnote{Bruce Berndt mentioned once that he was waiting for Dick's report for quite a while -- it had to be that the manuscript was buried on the bottom of that box!}} 
He worked behind the dinner table on the first floor next to the Ramanujan bust, standing in the corner (Figure~4), and to a personal copy machine, a precious gift from Paul Nevai. Dick loved to tell the well-known story of the creation of this bust, so I will not repeat it here; see also \cite{Askey90} for Askey's review dedicated to Ramanujan and Liz's excellent report on their trip to India, December 1987 -- January 1988 \cite{Tom}.
Part of Dick and Liz's house looked like a museum to me; I remember a thousand-year-old (or two thousand?) Chinese vase surrounded by numerous books on the second floor. Once, Dick picked up a photo album of great mathematicians \cite{Halmos}%
{\footnote{Many thanks to Shaun Cooper for identifying this book! (See Figure~5.)}}
from this collection and told me personal stories about them --
he knew most of them and I should have written those recollections down! 
It is very important for any researcher to be a part of the scientific community, to know ``what can be done and who can do it".
Members of the community, actually
working together, inspire and motivate each other with new ideas. 
All of that I gained from Professor Richard Askey.
Alexei Zhedanov wrote: ``I was born under the constellation of Askey" \cite[Contribution \#~85]{Liber}. For me Dick was a `center of gravity' in those astronomical terms.

\begin{figure}[h!]
\includegraphics[width=0.8\linewidth]{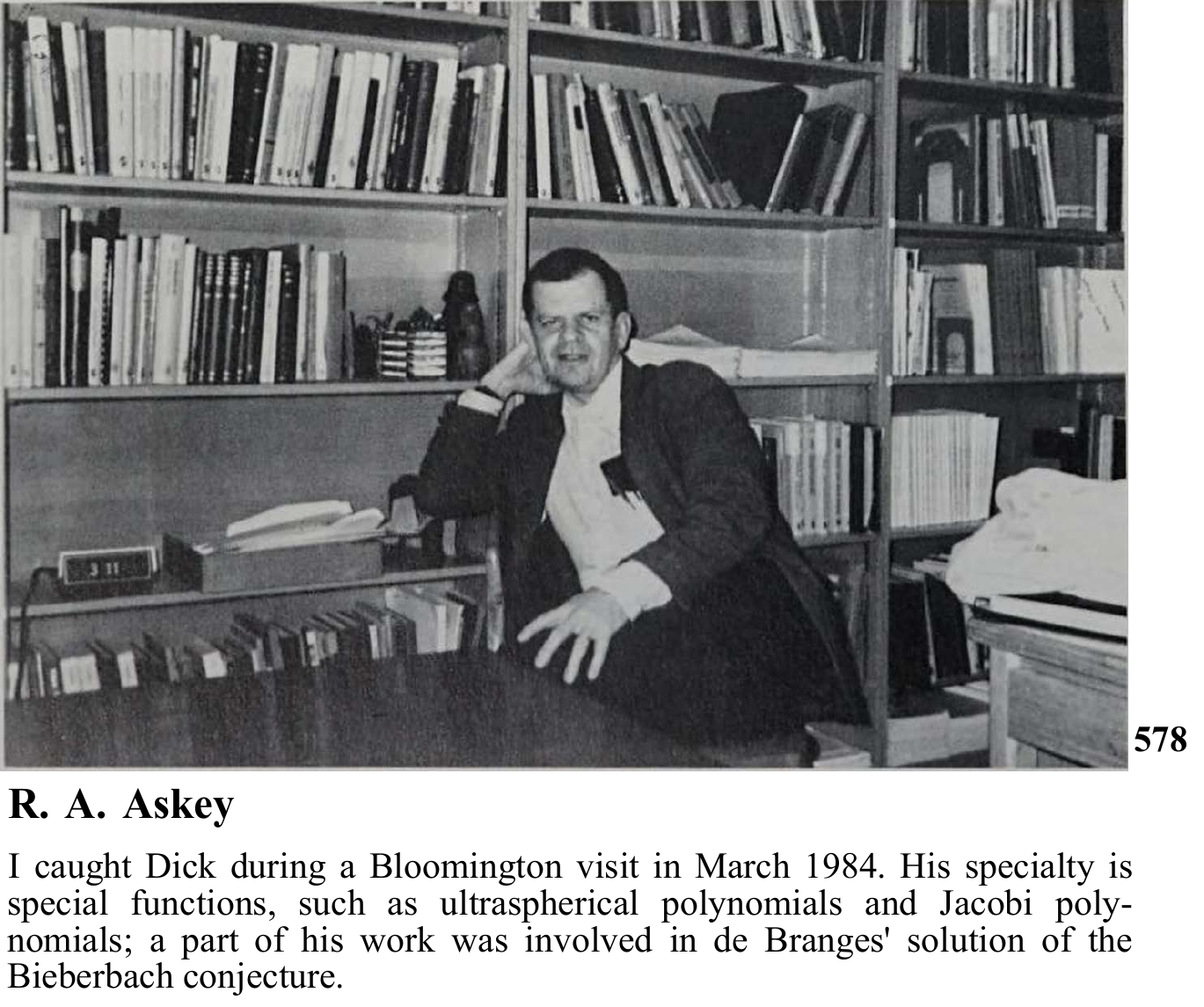} 
\caption{\cite[Photo~578]{Halmos} 
}
\end{figure}

I am very proud that our joint article with Mizan and Dick  \cite{ARS96}, which drew some attention recently, with only one known typo \cite{Centneretal, Coh:Sch} (thank you, Howard!), 
was initiated in Ottawa and had been completed in Dick's house among those treasures.

George Gasper kindly pick me up at the Askey's residence and drove me to Evanston, Illinois, 
before coming back to Moscow.

\section{In my apartment in Moscow and beyond\/}

During his second visit to Moscow, in 1991, Dick was invited to my apartment to meet the whole family.
He brought some flowers -- red roses. 
When Dick went to my children's room, which had a two-level bed for my daughter and son, he immediately said there was a similar bed at his house when he was a kid. Another recollection of his family was from the time of the Great Depression. As far as I remember, his father oversaw municipal transportation in a city during those hard times. When things went from bad to worse, instead of firing his employees, he said he would not fire anyone but would cut everyone's salary in half. Those were Dick's family traditions -- he was almost an American royalty for me and my students \cite[Contribution \#~71]{Liber}.

Dick's handling of necessary formalities was always quick.
One time when I had a problem obtaining an American visa in Moscow, Dick called directly to the embassy from his home. This was quite unusual in the former Soviet Union at that time.
He explained that I was on my way to Canada and asked them to expedite my case. Liz later told me in Madison that he didn't even record the name of the clerk at the embassy. The visa was issued the very next day! In the year 2000, when we had a similar problem with the former Soviet Union participants of The NATO Advanced Study Institute Special Functions 2000  {\cite{NATO2000, RMJM}}, following Dick's example, Joaqu\'{i}n Bustoz and I had to email a respectful request to the American embassy in Moscow, and all visas for our participants were also issued the next day! Such friendly relationships existed between the two countries! Unfortunately, it's just impossible nowadays when the state officials of two great countries did not talk to each other under the Biden administration!
During Dick's 
visit to Moscow, 
one morning, I found him on the second floor of the Steklov Institute, in the office of the associate director, receiving visitors, with Professor Naum Yakovlevich Vilenkin among them,
to the best of my recollection.
When I asked Askey how he was able to get into this prestigious office without any supervisor, Dick told me that he just asked the director's secretary, whose office is the next but one door in the corridor, 
and she opened this office for him and all his incoming visitors!
On another occasion, Dick had to make a copy of a document in Moscow where, at that time, it was difficult. 
What did he do? He asked a driver to stop at the American embassy and came back with 
the needed copy in 5 minutes! I was waiting in the car, but we could have gotten 
in big trouble at my work for this unexpected stop in that location. Once again, 
a direct solution if you have an American passport!
{\footnote{
According to Liz's diaries, they visited the U.S. Embassy before that in order to meet Paula Boyd in the Science and
Technology section \cite{Tom}.
}}

After that there were meetings/conversations/collaborations with Dick on different conferences 
all over the world.
In May of 1994, Luc Vinet organized a Workshop on Symmetries and Integrability of Difference Equations, Hotel Est\'{e}rel, Qu\'{e}bec, Canada \cite{Esterel}, as a ``mixture of special functions and nonlinear people". Two `patriarchs', Martin Kruskal and Richard Askey, were sitting in the hotel lobby, in the evenings, accepting numerous visitors/participants. 
Later, when I was originally hired as a senior lecturer by the Arizona State University 
in the late 1990s, my first permanent job in the United States, 
Dick immediately said that in the classroom I would always be a professor! 
Next year, when my NSF grant application had been approved, 
he called the analysis program director in order to tell him that he had made the right choice! 
There are many other ``real life stories" like that; for instance, Dick's
presentation of our talk in the International Congress of Mathematicians in Warsaw, when
Arnold Nikiforov ``couldn't attend" in 1983 (so did Professor Gel'fand!). 
\begin{figure}[h!]
\includegraphics[width=0.825\linewidth]{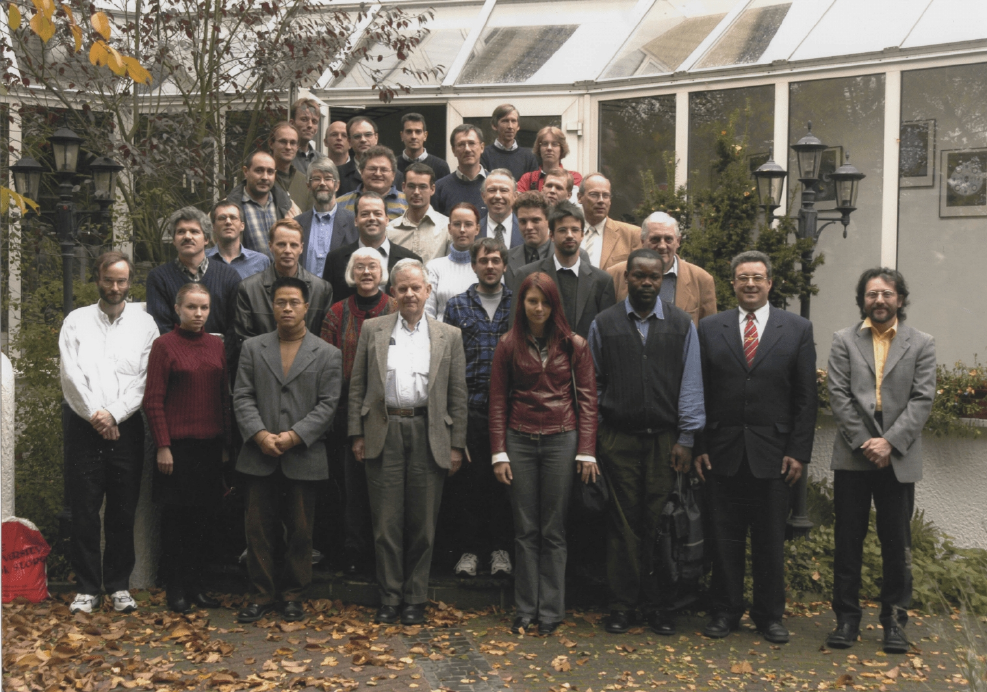} 
\caption{ Group photo at International Workshop on Special Functions,
Orthogonal Polynomials, Quantum Groups and Related Topics: dedicated to Dick Askey on his 70th Birthday,
Hotel Hochwiesm\"{u}hle, Bexbach, Saarland, Germany, October 18--22, 2003 \cite{bexbach}\/. }
\end{figure}
In Tempe, we had a chance to take Liz and Dick to the South Mountain park with a gorgeous view of the Phoenix metropolitan area from the elevation. Going to a Japanese restaurant after that, Dick went first in order to greet the host in Japanese. He liked this country very much (and talked about it a lot) after their trip to Japan in October--November 1987 \cite{Tom}. 
At home, Liz and Dick sometimes liked to listen to a relaxing classical music.
One day, Askey informed me that he will be interviewed by a radio station on mathematical education. 
This live event went well, Dick was very convincing -- but it was almost unrealistic to hear his familiar voice on the air from my house in Arizona!
In the United States, it was always exciting to discuss math and physics with Dick, 
such as just bumping into him during the Annual AMS--MAA Meetings.%
\footnote{By the way, I became a member of the AMS thanks to Dick in Russia -- he paid my membership fee until I came to the USA.\/}
In Washington D.~C., January 2000, for instance, he took me to 
the University of Wisconsin--Madison traditional alumni meeting 
with a great collection of only Wisconsin's cheeses; see also \cite{Askey00} for Dick's presentation. 
(There was a snow storm and, unlike Sicily, Dick's taxi had a minor accident on the way from the airport.)
In Boston, January 2013, I was fortunate to stop him for quite a while 
in the covered path between hotels and the convention center to tell him
that our first draft of an article on the minimum-uncertainty squeezed
states, one of the fundamental concepts in quantum physics, had been rejected by
Journal of Physics~B in just 24 hours! Obviously, he found some words to
encourage us and an extended revision was later published \cite%
{MinUnSqSt2013}. This pretty long conversation had ended when Dick mentioned that he was late to meet Liz.
And this is just one episode of Dick's friendly support of
our research; he always has been an uncompromising judge of novelty and
quality of publications in many areas of analysis and later in mathematical
education. 

At the beginning of his well-received talk at Arizona State University, more then a decade ago \cite{Askey09}, 
Professor Askey said: ``we cannot lie to our students but we cannot tell them the whole truth".
After the talk, he shared with our faculty some details on `his' presentation in Warsaw, 1983.
In the context of `rigorous proofs', he mentioned once that ``if it doesn't work formally, it doesn't work".
It is so sad that there is no way to present him a copy of our recent article \cite{PhUs2022}, 
which has attracted certain attention in view of the 100th anniversary of the birth of quantum mechanics.
Dick called me a physicist and I have always been wondering what else should I do in order to be recognized by Askey as a mathematician? Now; see also Figure~2; one may say that he was absolutely right!
According to an {\textsl{AI Overview\/}}, generated on January 29, 2025
\footnote{Please google: dick askey, special functions, applications, quantum physics\ .}
:
``Dick Askey is a prominent mathematician well-known for his significant contributions to the field of `special functions', particularly their applications in various areas of mathematics and physics, including quantum mechanics; where special functions often appear as solutions to key equations describing quantum systems
...\ 
Askey's contribution:
Dick Askey is recognized for his extensive research on the theory of special functions, particularly on developing new relationships and identities between different special functions, which has greatly enhanced their applicability in various fields, including quantum physics\ ..." (One may consider \cite{As:Su1, As:Su2}, as an example? -- SKS.) The article \cite{PhUs2022} is helping the readers to study quantum physics from a historical perspective.
Last but not least, we remember the International Workshop on Special Functions, Orthogonal Polynomials, Quantum Groups, and Related Topics, dedicated to Askey's 70th birthday. It was organized by Andreas Ruffing in Bexbach, Germany \cite{And:Gasp:Sus, bexbach} (Figure~6), featuring incredible fireworks in Dick's honor; see also \cite{Ey:Ruff:Sus} and \cite{SusFS}.%
{\footnote{From Ruffing's email to Suslov on 11th January 2025:
``It is a pity that you could not be at the castle of Otocec, Slovenia, 2003, at the Autumn School of AbiTUMath where I took Dick and Liz Askey, having started from Munich: namely these
were the days before the Bexbach conference. There was a concert in Otocec, at the Autumn School, of a famous Slovenian string quartet.
The concert had been organized and the contact to the quartet had been provided by our friend Marko Robnik.
We had ordered the quartet in honor of the guests, and very particular in honor of Dick.
Dick did like the quartet so much, and it was, I think, the only time I observed that he was close to tears.
He said later `Special thanks to the one who has invited this impressive quartet, an extraordinary
idea.'
This evening was some kind of Advent for Bexbach 2003."}}

Those feelings of freedom, support, expertise, and joy of life, coming from Dick, are my major personal recollections of precious moments with the distinguished mathematician, Professor Richard Askey!

In the end, everyone faces the question: ``Do you love science or yourself in science?" In other words, are you pursuing scientific research with genuine sincerity, passion for the subject itself, and not just for personal gain, career promotion or recognition? 
Dick was a true mathematician of the first kind!
%

%
%

\smallskip 
\noindent \textbf{Acknowledgments.\/} I am grateful to Ruben Abagyan, Kamal Barley, Howard Cohl, Ricardo Cordero-Soto, Shaun Cooper, Mourad Ismail, Tom Koornwinder, Sergey Kryuchkov, Nathan Lanfear, Peter Paule, Andreas Ruffing, Dennis Stanton, Eugene Stepanov, and Lilya Suslova for their valuable help.

\bigskip 
{\small {{{\scshape{Orcid:\/}} {{\url{https://orcid.org/0000-0001-8169-0987}}%
}}}}


\begin{thebibliography}{99}
\bibitem{An:As:RoySF} G.~E.~Andrews, R.~Askey, and R.~Roy, \textsl{Special
functions\/}, Encyclopedia of Mathematics and its Applications, vol.~71,
Cambridge University Press, Cambridge, 1999. {%
\url{https://doi.org/10.1017/CBO9781107325937}}

\bibitem{And:Gasp:Sus} G.~E.~Andrews, G.~Gasper, and S.~K.~Suslov, \emph{Preface\/}. The Ramanujan Journal,
Askey Special Issues (G.~E.~Andrews, G.~Gasper, and S.~K.~Suslov, Coordinating Editors), \textbf{13} (2007), 5--6. 
{\url{https://link.springer.com/article/10.1007/s11139-006-0239-z}}

\bibitem{Askey75} R.~Askey, \textsl{Orthogonal polynomials and special
functions\/}, Society for Industrial and Applied Mathematics, Philadelphia,
Pa., 1975.

\bibitem{Askey85} R.~Askey, \emph{Continuous Hahn polynomials\/},  Journal
of Physics A: Mathematical and General \textbf{18} (1985)~\#16,
L1017--L1019. {%
\url{https://iopscience.iop.org/article/10.1088/0305-4470/18/16/004}}

\bibitem{Askey90} R.~ Askey, \emph{Ramanujan and hypergeometric and basic hypergeometric series\/} [in Russian], translated from English with remark by N.~M.~Atakishiyev and S.~K.~Suslov, Uspekhi Mat. Nauk \textbf{45} (1990)~\#1, 33--76; Russian Mathematical Surveys \textbf{45} (1990)~\#1, 37--86, The British Library Board and The London Mathematical Society.
    {\url{https://iopscience.iop.org/article/10.1070/RM1990v045n01ABEH002325}}
    

\bibitem{Askey00} R.~A.~Askey, \emph{School mathematics CDs from Singapore\/},  MAA Special Presentation,  Washington~D.~C., January~22nd, 2000. 
    {\url{https://jointmathematicsmeetings.org/meetings/national/jmm/2026_progfull.html}}

\bibitem{Askey09} R.~Askey, \emph{Algebra, geometry and trigonometry\/}, Colloquium: Research Innovations in Mathematics and Science Education in the School of Mathematics and Statistics at Arizona State University, February~17th, 2009. \\
     \textsl{Abstract\/}: {Serious trigonometry started with Ptolemy's theorem. Three or four proofs of this theorem will be given, using algebra and trigonometry. There are four standard concurrence theorems for a triangle. The harder two are the medians and the altitudes. Both of these will be considered\/}.

\bibitem{ARS96} R.~A.~Askey, M.~Rahman, S.~K.~Suslov, \emph{On a general $q$%
-Fourier transformation with nonsymmetric kernels\/}, Journal of
Computational and Applied Mathematics \textbf{68} (1996)~\#1--2, 25--55.  {%
\url{https://doi.org/10.1016/0377-0427(95)00259-6}} 

\bibitem{As:Su1} R.~Askey and S.~K.~Suslov, \emph{The $q$-harmonic oscillator and an analogue of the Charlier polynomials\/}, Journal of Physics A: Mathematical and General \textbf{26} (1993)~\#15, L693--L698.
    {\url{https://iopscience.iop.org/article/10.1088/0305-4470/26/15/014}}

\bibitem{As:Su2} R.~Askey and S.~K.~Suslov, \emph{The $q$-harmonic oscillator and the Al-Salam and Carlitz polynomials\/}, Letters in Mathematical Physics \textbf{29} (1993), 123--132.
    {\url{https://link.springer.com/article/10.1007/BF00749728}}

\bibitem{AW85} R.~Askey and J.~Wilson, \emph{Some basic hypergeometric
orthogonal polynomials that generalize Jacobi polynomials\/}, Memoirs of the
American Mathematical Society \textbf{54},(1985)~\# 319, iv+55~pp. {%
\url{https://www.ams.org/books/memo/0319/}}

\bibitem{At:SusCH} N.~M.~Atakishiyev and S.~K.~Suslov, \emph{The Hahn and
Meixner polynomials of an imaginary argument and some of their applications\/%
}, Journal of Physics A:  Mathematical and General \textbf{18} (1985)~\#10, 1583--1596.  {%
\url{https://iopscience.iop.org/article/10.1088/0305-4470/18/10/014}} 

\bibitem{At:SusAW} N.~M.~Atakishiyev and S.~K.~Suslov, \emph{On the
Askey--Wilson polynomials\/}, Constructive Approximation \textbf{8} (1992),
363--369.  {\url{https://doi.org/10.1007/BF01279025}} 


\bibitem{PhUs2022} K.~Barley, J.~Vega-Guzman, A.~Ruffing, and S.~K.~Suslov, 
\emph{Discovery of the relativistic Schr\"{o}dinger equation\/},
Physics-Uspekhi \textbf{65} (2022)~\#1, 90--103. From the history of
physics. (This article has been selected by The Institute of Physics=IOP as
one(\#3) of the top 5 articles in a January 2023 collection entitled ``Dive
deeper into mathematics and computational research") One of the most read
(\#1-\#4) articles in Physics-Uspekhi since April 2022 until May 2023.  {%
\url{https://iopscience.iop.org/article/10.3367/UFNe.2021.06.039000}} 

\bibitem{bexbach} The Bexbach meeting, Germany, October 2003. {\url{https://math.la.asu.edu/~suslov/bexbach/index.html}}

\bibitem{NATO2000} J.~Bustoz, M.~E.~H.~Ismail, and S.~K.~Suslov, Editors, \textsl{Special Functions 2000: Current Perspective and Future Directions\/} (NATO Science Series II: Mathematics, Physics and Chemistry, vol. 30, Kluwer Academic Publishers/Springer-Verlag, Dordrecht-Boston-London, 2001.
    {\url{https://link.springer.com/book/10.1007/978-94-010-0818-1}}

\bibitem{RMJM} J.~Bustoz, M.~E.~H.~Ismail, and S.~K.~Suslov, Guest Editors, Rocky Mountain Journal of Mathematics, Special Issue on Special Functions: Proceedings of the NATO Advanced Institute Special Functions 2000, vol.~32, no.~2, 2002, pp. 1--936.
      
\bibitem{Centneretal} R.~Centner, H.~S.~Cohl, and R.~S.~Costas-Santos,
    \emph{The Askey--Rahman--Suslov nonsymmetric Poisson kernel for the Askey--Wilson polynomials and its special values\/}, to appear.
   
\bibitem{Chud} The Chudnovsky brothers, Wikipedia,
{\url{https://en.wikipedia.org/wiki/Chudnovsky_brothers}}

\bibitem{Liber} H.~S.~Cohl and M.~E.~H.~Ismail, Editors, \textsl{Liber Amicorum, Richard ``Dick'' Allen Askey -- a Friendship Book -- from Dick's colleagues and friends (September 15, 2019)\/},
 Celebratio Mathematica, 2022.
{%
\url{https://celebratio.org/Askey_RA/article/1009/}%
}

\bibitem{RAleg} H.~S.~Cohl, M.~E.~H.~Ismail, and H.-H.~Wu, \emph{The Legacy
of Dick Askey (1933--2019)\/}, Notices of the American Mathematical Society 
\textbf{69} (2022)~\#1, 59--75. {%
\url{https://www.ams.org/notices/202201/rnoti-p59.pdf}}

\bibitem{Coh:Sch} H.~S.~Cohl and M.~J.~Schlosser,
    \emph{Quadratic and linear transformation formulas for nonterminating basic hypergeometric series 
    by evaluations of Askey--Wilson polynomials\/}, to appear.

\bibitem{Ey:Ruff:Sus} K.~Ey, A.~Ruffing, and S.~K.~Suslov, 
\emph{Method of separation of the variables for basic analogs of 
equations of mathematical physics\/}, The Ramanujan Journal, 
Askey Special Issues (G.~ E.~Andrews, G.~Gasper, and S.~K.~Suslov, Coordinating Editors), 
\textbf{13} (2007) \#1--3, 407--447.
{\url{https://link.springer.com/article/10.1007/s11139-006-0260-2}}

\bibitem{Faddeev} L.~D.~Faddeev, {\textsl{Autobiography of Ludwig Faddeev\/}}, The Shaw Prize in Mathematical Sciences, 9~September 2008, Hong Kong.
     {\url{https://www.shawprize.org/autobiography/ludwig-faddeev/}}; \;
    {\textsl{Ludwig D. Faddeev. Obituary by Nikolai Reshetikhin, Michael Semenov-Tian-Shansky, and Leon Takhtajan\/}}, IAMP News Bulletin, October 2017, pp.~29--45.   
    {\url{https://www.math.stonybrook.edu/~leontak/L.D.%20Faddeev-IAMP%20News%20Bulletin%20October%202017.pdf}}

\bibitem{GelGr} I.~M.~Gelfand and M.~I.~Graev, \emph{$GG$-functions and
their relations to ${\mathcal{A}}$-hypergeometric functions\/}, 
arXiv:math/9905134v1 [math.AG] 20 May 1999.
{\url{https://arxiv.org/pdf/math.AG/9905134}}

\bibitem{Halmos} P.~R.~Halmos, \textsl{I Have a Photographic Memory\/}, 
    American Mathematical Society, Providence, Rhod Island, 1987.
    
\bibitem{Esterel}  The Hotel Est\'{e}rel, Qu\'{e}bec, Canada.
{\url{https://www.hotelsone.com/esterel-hotels-ca/esterel-resort.html}}; \;
{\url{https://www.esterel.com/en/business-meetings-and-events/}} 


\bibitem{Tom} T.~Koornwinder, \emph{Dick and Liz Askey's world trip in fall of 1987\/},
    Lecture by Tom Koornwinder, Diary texts read by Suzanne Askey, 16th International Symposium on Orthogonal Polynomials, Special Functions and Applications, online, June 13--17, 2022.
    {\url{https://staff.fnwi.uva.nl/t.h.koornwinder/art/sheets/2022_OPSFA.pdf}}; \\
    Liz's edited originals: \\
    {\url{https://staff.fnwi.uva.nl/t.h.koornwinder/specfun/#DickAskey}};\; \\
    {\url{https://staff.fnwi.uva.nl/t.h.koornwinder/specfun/AskeyDiary_USSR1987.pdf}};\; \\
    {\url{https://staff.fnwi.uva.nl/t.h.koornwinder/specfun/AskeyDiary_Japan1987.pdf}};\; \\
    {\url{https://staff.fnwi.uva.nl/t.h.koornwinder/specfun/AskeyDiary_Australia1987.pdf}}; \; \\
    {\url{https://staff.fnwi.uva.nl/t.h.koornwinder/specfun/AskeyDiary_India1987.pdf}}
    
\newpage

\bibitem{MinUnSqSt2013} S.~I.~Kryuchkov, S.~K.~Suslov, and
J.~M.~Vega-Guzman, \emph{The minimum-uncertainty squeezed states for atoms
and photons in a cavity\/}, Journal of Physics B: Atomic, Molecular and
Optical Physics \textbf{46} (2013)~\#10, 104007 (15~pp). (IOP=Institute Of
Physics SELECT and HIGHLIGHT for 2013).  {%
\url{https://iopscience.iop.org/article/10.1088/0953-4075/46/10/104007}} 
    

\bibitem{Kurchatov} The Kurchatov Institute, Moscow. {\url{http://nrcki.ru/}} ;
see also {\url{https://commons.wikimedia.org/wiki/File:%D0%9F%D0%B0%D0%BC%D1%8F%D1%82%D0%BD%D0%B8%D0%BA_%D0%98._%D0%92._%D0%9A%D1%83%D1%80%D1%87%D0%B0%D1%82%D0%BE%D0%B2%D1%83_2021.jpg}}

\bibitem{Lan:Lif} L.~D.~Landau and E.~M.~Lifshitz, \textsl{Mechanics\/}, 3rd Edition,
Course in Theoretical Physics, vol.~1, Butterworth-Heinemann, 2005.

\bibitem{FizTex} The Moscow Institute of Physics and Technology. {\url{https://eng.mipt.ru/why-mipt/}}

\bibitem{NSU} A.~F.~Nikiforov, S.~K.~Suslov, and V.~B.~Uvarov, \textsl{%
Classical Orthogonal Polynomials of a Discrete Variable\/} [in Russian],
Nauka, Moscow, 1985; English translation in Springer Series in Computational
Physics, Springer-Verlag, 1991.

\bibitem{Petergof} The Peterhof Palace, St.~Petersburg. 
{\url{https://en.wikipedia.org/wiki/Peterhof_Palace}}

\bibitem{Rahman} Mizan Rahman, Wikipedia, 
    {\url{https://en.m.wikipedia.org/wiki/Mizan_Rahman}}

\bibitem{AtomicBomb} R.~Rhodes, \textsl{The Making of the Atomic Bomb\/}, 1st Edition, Simon \& Schuster, 1987. 
{\url{https://www.abebooks.com/servlet/BookDetailsPL?bi=32091573561}}

\bibitem{Smor} Ya.~A.~Smorodinskii, \textsl{Selected Works\/}, 2nd Edition, Classics of Science, Editorial URSS, Moscow,  2006 [in Russian].

\bibitem{SteklovMoscow} The Steklov Mathematical Institute of Russian Academy of Sciences, Moscow. {\url{https://www.mi-ras.ru/?l=1}}

\bibitem{SteklovStP} Saint Petersburg Department of V.~A.~Steklov Institute of Mathematics of the Russian Academy of Sciences.
    {\url{https://www.pdmi.ras.ru/eng/institut/institut.php}} \;
    {\url{https://www.pdmi.ras.ru/pdmi/en/EIMI}}

\bibitem{SusFS} S.~K.~Suslov, \textsl{An Introduction to Basic Fourier Series\/}, 
    Kluwer Series ``Developments in Mathematics", vol.~9, Kluwer Academic Publishers/Springer-Verlag, 
    Dordrecht-Boston-London, 2003.
    {\url{https://math.la.asu.edu/~suslov/bfs/}}; see also
    {\url{https://math.la.asu.edu/%7Esuslov/bfs/reiher.pdf}}\ and
    {\url{https://math.la.asu.edu/%7Esuslov/bfs/bfserrata.pdf}}
    
\bibitem{Trinity} Trinity United Methodist Church, Madison, Wisconsin.
    {\url{https://www.tumcmadison.com/}}   
    
\bibitem{Wilson91}  J.~A.~Wilson, \emph{Orthogonal functions from Gram determinants\/},
    SIAM Journal on Mathematical Analysis \textbf{22} (1991)\#~4, 10.1137/0522074.
    {\url{https://epubs.siam.org/doi/10.1137/0522074}}   

\end{thebibliography}
\end{document}